\documentclass[envcountsame,envcountsect]{svmult}
\usepackage{mathdots}
\usepackage{graphicx}
\usepackage{bm}
\usepackage{mathtools}
\usepackage{mathrsfs}
\usepackage{type1cm}
\usepackage{latexsym,enumitem}
\usepackage{dsfont}
\usepackage{bbm}
\usepackage{soul}
\usepackage{epic,tikz}
\usepackage{xcolor}

\usepackage{latexsym}
\usepackage{dsfont}
\usepackage{bbm}
\usepackage{amssymb}
\usepackage{amsmath}

\renewcommand {\epsilon}{\varepsilon}
\renewcommand {\le}{\leqslant}
\renewcommand {\ge}{\geqslant}
\renewcommand {\le}{\leqslant}
\renewcommand {\geq}{\geqslant}

\RequirePackage[colorlinks,citecolor=blue,urlcolor=blue]{hyperref}

\newtheorem{Theorem}{Theorem}[section]

\newtheorem{Cor}[Theorem]{Corollary}

\newtheorem{Rem}[Theorem]{Remark}

\newtheorem{Def}[Theorem]{Definition}

\newcommand{\ssup}[1] {{\scriptscriptstyle{({#1}})}}

\def\cB{\mathscr{B}}

\def\cD{\mathscr{D}}

\def\cI{\mathscr{I}}

\def\cM{\mathscr{M}}

\def\cP{\mathscr{P}}

\def\cS{\mathscr{S}}

\def\fS{\mathfrak{S}}

\def\Erw{\mathbb{E}} 

\def\G{\mathbb{G}}

\def\N{\mathbb{N}}
  
\def\Prob{\mathbb{P}} 
\def\Q{\mathbb{Q}}
\def\R{\mathbb{R}}

\def\X{\mathbb{X}}

\def\Z{\mathbb{Z}}

\def\bX{\textit{\bfseries X}}

\def\eps{\varepsilon}
\def\vph{\varphi}

\def\1{\mathbf{1}}
\def\3{{\ss}}

\def\weakly{\stackrel{w}{\to}}
\def\RA{\Rightarrow}

\def\wh{\widehat}

\def\ovl{\overline}

\def\sign{\textsl{sign}}

\def\dd{\mathrm{d}}

\allowdisplaybreaks[4] 

\begin{document}

\title*{\Large ON NULL-HOMOLOGY AND STATIONARY SEQUENCES}
\titlerunning{\large\sc Null-homology and stationary sequences}
\author{\large\sc By Gerold Alsmeyer and Chiranjib Mukherjee}
\institute{Inst.~Math.~Stochastics, Department
of Mathematics and Computer Science, University of M\"unster,
Orl\'eans-Ring 10, D-48149 M\"unster, Germany.\at
\email{gerolda@math.uni-muenster.de, chiranjib.mukherjee@uni-muenster.de}\at
Funded by the Deutsche Forschungsgemeinschaft (DFG) under Germany's Excellence Strategy EXC 2044--390685587, Mathematics M\"unster: Dynamics--Geometry--Structure.}

\maketitle

\abstract{
The concept of homology, originally developed as a useful tool in algebraic topology, has by now become pervasive in quite different branches of mathematics. The notion particularly appears quite naturally in ergodic theory in the study of measure-preserving transformations arising from various group actions or, equivalently, the study of stationary sequences when adopting a probabilistic perspective as in this paper. Our purpose is to give a new and relatively short proof of the coboundary theorem due to Schmidt \cite{Schmidt:77} which provides a sharp criterion that determines (and rules out) when two stationary processes belong to the same \emph{null-homology equivalence class}. 
We also discuss various aspects of null-homology within the class of Markov random walks, compare null-homology with a formally stronger notion which we call {\it strict-sense null-homology}. Finally, we also discuss some concrete cases where the notion of null-homology turns up in a relevant manner. }

\bigskip



\section{Introduction and motivation}\label{sec:intro}

The notion homology arises in various branches of mathematics. It was originally developed in algebraic topology in order to associate a sequence of algebraic objects. A typical fundamental question is the following: When does an $n$-cycle of a (simplical) complex form the boundary of an $(n+1)$-chain,  or equivalently, when is its fundamental class  a boundary for the singular homology? If such a requirement is fulfilled, the cycle is said to be \emph{homologous to $0$} or \emph{null-homologous}. In ergodic theory, the notion appears in the context of measure-preserving transformations arising from natural group actions on a complete and separable metric space. A prominent result in this context, due to K. Schmidt \cite[Theorem 11.8]{Schmidt:77} and sometimes called the \emph{coboundary theorem}, states that a random walk or \emph{cocycle} generated by a measure-preserving transformation of $\G=\R^{d}$ or a closed subgroup of it is tight iff it is null-homologous (and then called a \emph{coboundary}). The result was extended  to locally compact second countable Abelian groups $\G$ with no compact subgroup in \cite{MooreSchm:80}, and generalized to Polish topological groups $\G$ in \cite{AarWeiss:00}. In probabilistic language, the coboundary theorem provides a characterization of tightness for a random walk with stationary increments. It will be stated below both in the ergodic theoretic and the probabilistic framework. A similar framework was also used by Bradley who, by adapting Schmidt's arguments, proved the result for random walks with  nonstationary $\G$-valued increments; first for $\G=\R$ in \cite{Bradley:95}, then for $\G$ a separable Banach space in \cite{Bradley:97} and for products of upper triangular random matrices in \cite{Bradley:96}.

\vspace{.1cm}
In order to motivate the present work, let us note that a random walk with stationary increments can be viewed as a particular instance of a so-called \emph{Markov random walk} or \emph{Markov-additive process} with stationary recurrent driving chain. A precise definition will be given shortly. Our purpose is to give a new and relatively short proof of the aforementioned coboundary theorem (see Theorem \ref{thm-coboundary}), discuss various aspects of null-homology within the class of Markov random walks (in Section \ref{sec-MRW}), compare null-homology with a formally stronger notion we call {\it strict-sense null-homology} (cf. Theorem \ref{thm-main 1})
and provide an extension of the coboundary theorem to the case when the driving chain is null-recurrent (cf. Theorem \ref{thm-main 2}. Finally, we also provide further background information and applications (cf. Section \ref{sec-appli}) so as to sustain the relevance of these results in probability theory where the notion of null-homology has remained relatively unnoticed, at least until its appearance in recent progress made in the {\it strong coupling limit of the Polaron path measures} \cite{MV:18} (see Section \ref{sec-polaron}) and in stochastic homogenization (see Section \ref{sec-homogenization}) under the name increment-stationarity, see e.g.~Blanc, Le Bris \& Lions \cite{BlancLeBLions:07}, Gloria \& Otto \cite{GloriaOtto:11} and especially Gloria \cite{Gloria:14}.

\section{The coboundary theorem, revisited}\label{sec-coboundary}

We start with the ergodic-theoretic framework, where the coboundary theorem is usually cast in terms of a measure-preserving transformation and a cocycle. 
In this framework we give a new and short proof of the above result with the help of commuting maps and Schauder's fixed-point theorem.

\vspace{.1cm}
We proceed with a definition of null-homology in terms of probability measures rather than random variables. Let $\Omega=(\R^{d})^{\otimes\Z}$ be the space of doubly infinite sequences ${\bf x}=(x_{n})_{n\in\Z}$ endowed with the Borel $\sigma$-field and $T:\Omega\to\Omega$ the (left) shift operator on $\Omega$, viz.
$$ {\bf x}\ =\ (\ldots,x_{-1},x_{0},x_{1},\ldots)\ \mapsto\ (\ldots,x_{0},x_{1},x_{2},\ldots). $$
The coordinate mappings on $\Omega$ are denoted $X_{n}$ for $n\in\Z$, and we let $S_{n}$ be the mapping ${\bf x}\mapsto s_{n}$ on $\Omega$ for $n\in\Z$, where
\begin{align}\label{def-nullhom-measure1}
s_{n}\ =\ \begin{cases} \hfill x_{1}+\ldots+x_{n}&\text{if }n\ge 1,\\
\hfill 0&\text{if }n=0,\\
-(x_{-n+1}+\ldots+x_{0})&\text{if }n\le -1.
\end{cases}
\end{align}

\vspace{.1cm}
Next, let $\cM(\Omega)$ denote the locally convex vector space of finite signed measures on $\Omega$ endowed with the topology of weak convergence and $\cM_{T}(\Omega)$
its subset of $T$-invariant 
probability measures. Defining the map $D:\Omega\to\Omega$ by 
$$ {\bf x}\ \mapsto\ T{\bf x}-{\bf x}\ =\ (\ldots,x_{0}-x_{-1},x_{1}-x_{0},x_{2}-x_{1},\ldots), $$
we obviously have that $\Prob\in\cM_{T}(\Omega)$ implies $\Prob D^{-1}:=\Prob(D\in\cdot)\in\cM_{T}(\Omega)$. Null-homology for elements of $\cM_T(\Omega)$ can now be defined as follows. 
\begin{Def}\label{def-nullhom-measure2}
A 
$T$-invariant probability measure $\Prob\in\cM_T(\Omega)$ is called null-homologous if $\Prob=\Q D^{-1}$ for some $\Q\in\cM_T(\Omega)$.
\end{Def}


\begin{Theorem}\label{thm-coboundary}
Given any $\Prob\in\cM_{T}(\Omega)$, the following assertions are equivalent:
\begin{itemize}\itemsep2pt
\item[(a)] $\{\Prob S_{n}^{-1}:n\ge 0\}$ is tight.
\item[(b)] $\Prob$ is null-homologous.
\end{itemize}
\end{Theorem}

\begin{proof}
It suffices to show that (a) implies (b) for which we consider the bivariate mappings
\begin{gather*}
\Lambda_{k}:\Omega\to\Omega\times\Omega,\\
{\bf x}\ \mapsto\ \big(T^{k}{\bf x},{\bf x}+\ldots+T^{k-1}{\bf x}\big)\ =\ \big((x_{n+k})_{n\in\Z},(x_{n}+\ldots+x_{n+k-1})_{n\in\Z}\big)
\end{gather*}
for $k\in\N$ and point out that (b) entails the tightness of the family
$$ \cP\ =\ \{\Prob\Lambda_{k}^{-1}:k\in\N\}. $$
We can lift the shift $T$ as well as the projections $X_{n}$ in a canonical way to mappings on $\Omega\times\Omega$ and, by slight abuse of notation, may call these mappings again $T$ and $X_{n}$. The projections on the $y$-components, namely $({\bf x},{\bf y})\mapsto y_{n}$ if ${\bf y}=(y_{k})_{k\in\Z}$, are denoted $Y_{n}$ for $n\in\Z$. Then the $T$-invariance of $\Prob$ implies the very same for the elements of $\cP$. 

Now let $\cD$ be the closed convex hull of all weak limit points of $\cP$ which forms a compact convex subset of $\cM(\Omega\times\Omega)$. Consider the map
$$ S:\Omega\times\Omega\to\Omega\times\Omega,\quad \big({\bf x},{\bf y}\big)\ \mapsto\ \big(T{\bf x},{\bf x}+{\bf y}\big) $$
which is linear, continuous, commutes with $T$, i.e.~$S\circ T=T\circ S$, and satisfies further
$S\circ\Lambda_{n}=\Lambda_{n+1}$, thus $\Gamma_{n}S^{-1}=\Gamma_{n+1}$ for all $n\in\N$, where $\Gamma_{n}:=\Prob\Lambda_{n}^{-1}$. Then the last property entails that the set $\cD$ is $S$-invariant which in turn, by invoking Schauder's fixed point theorem,  allows us to conclude that $S$ has a fixed point, say $\Gamma$, in $\cD$. This means that $\Gamma S^{-1}=\Gamma$ or, equivalently, that $\Gamma$ is $S$-invariant.

Finally, by considering the map
$$ G=(X_{0},Y_{0}):\Omega\times\Omega\to\R^{d}\times\R^{d},\quad \big({\bf x},{\bf y}\big)\ \mapsto\  (x_{0},y_{0}) $$
we have that $(X_{n}',Y_{n}'):=G\circ S^{n}=(X_{n},Y_{n}\circ S^{n})$, $n\ge 0$, is stationary under $\Gamma$ and satisfies:
\begin{itemize}\itemsep2pt
\item[(1)] $(X_{n}')_{n\ge 0}=(X_{n})_{n\ge 0}$ and has law $\Prob$ under $\Gamma$, because this is the case under any element of $\cD$.
\item[(2)] $Y_{n+1}'=Y_{0}+X_{0}+\ldots+X_{n}=Y_{n}'+X_{n}'$, thus
$$ X_{n}\ =\ X_{n}'\ =\ Y_{n+1}'-Y_{n}' $$
for all $n\ge 0$.
\end{itemize}
Since $(Y_{n}')_{n\ge 0}$ is stationary under $\Gamma$, (a) follows.
\end{proof}

\begin{Rem}\rm
As mentioned earlier, the coboundary theorem was earlier proved by Schmidt \cite[Theorem 11.8]{Schmidt:77} using ergodic decomposition of ``skew-products". The technique there relied on showing that a cocycle is a coboundary if and only if its skew product decomposes into finite measure-preserving transformations. The above proof of Theorem \ref{thm-coboundary} is quite different (and a bit shorter) compared to the arguments in \cite[Ch. 11]{Schmidt:77}.
 \end{Rem}

\section{Markov random walks, null-homology and strict-sense null-homology}\label{sec-strict}
Let $\bX=(X_{n})_{n\in \Z}$ denote a doubly infinite stationary sequence of random variables which are defined on a probability space with underlying probability measure $\Prob$ and take values in a complete separable metric space $\cS$. Stationarity means that, for all $n\in\N$ and $m\in \Z$,
\begin{equation*}
\Prob\big((X_{1},\dots, X_{n}) \in \cdot\big)\ =\ \Prob\big((X_{m+1},\dots, X_{m+n})\in\cdot\big).
\end{equation*}
In other words, the joint law of $(X_{1},\dots, X_{n})$ for any $n$ coincides with the law of any of its ``shifts"  under the action of the additive group $\Z$ on the space of doubly-infinite sequences $\mathscr S^{\times\Z}$. 
Equivalently, putting $\bX_{n}:=(X_{n+k})_{k\in\Z}$ for $n\in\Z$, the law of $\bX_{n}$ is the same for each $n\in\Z$.
The notion of \emph{homology} now arises naturally from this group action, given measurable functions $F, G\colon \mathscr S^{\times\Z}\to \R^{d}$. Following Lalley \cite[p.~197]{Lalley:86}, we say that $F$ is \emph{homologous} to $G$ (with respect to $\bX$ and $\Prob$) and write $F\sim G$ if there exists a function $\xi: \mathscr S^\Z\to \R^{d}$ such that
\begin{equation}\label{eq:F,G homologous}
F(\bX_{1}) - G(\bX_{1})\ =\ \xi (\bX_{1})- \xi(\bX_{0})\quad\Prob\text{-a.s.}
\end{equation}
Then $\sim$ is an equivalence relation, and $F$ is called \emph{null-homologous} if $F\sim 0$, thus
\begin{equation}\label{eq:F null-homologous}
F(\bX_{1})\ =\ \xi (\bX_{1})- \xi(\bX_{0})\quad\Prob\text{-a.s.}
\end{equation}
Observe that, given any stationary sequence $\bX$ and null-homologous function $F$, the process $(F(\bX_{n}))_{n\in\Z}$ is not only stationary as well (and thus also tight), but in fact the incremental sequence of another stationary process, viz. $(\xi(\bX_{n}))_{n\in\Z}$. In view of this in fact equivalent definition of null-homology, the coboundary theorem answers the natural question which stationary processes are of that ``incremental'' type and therefore allowing a representation with respect to a null-homologous function.

\subsection{Markov random walks}\label{sec-MRW}


In \cite{Lalley:86}, Lalley considered random walks with increments from a fairly general class of \emph{integrable} stationary sequences. As a main result, he proved a Blackwell-type renewal theorem for which it was necessary to rule out a certain ``lattice-type" behavior which to define requires the notion of null-homology. In the following, we give a brief introduction of this notion within the more general framework of \emph{Markov random walks}.

\vspace{.1cm}
Let $(\cS,\fS)$ be a measurable space and $\cB(\R^{d})$ the Borel $\sigma$-field on $\R^{d}$. A sequence $(M_{n},X_{n})_{n\ge 0}$ taking values in $(\cS\times\R^{d},\fS\otimes\cB(\R^{d}))$ is called \emph{Markov-modulated} if
\begin{itemize}\itemsep2pt
\item $(M_{n})_{n\ge 0}$ forms a temporally homogeneous Markov chain,
\item $X_{0},X_{1},\ldots$ are conditionally independent given $(M_{n})_{n\ge 0}$, and
\item the conditional laws of $X_{0}$ and $X_{n}$ for $n\ge 1$ given that $M_{0}=s_{0},M_{1}=s_{1},\ldots$ are given by $P_{0}(s_{0},\cdot)$ and $P((s_{n-1},s_{n}),\cdot)$, respectively, for suitable stochastic kernels $P_{0}$ and $P$.
\end{itemize} 
It follows that $(M_{n},X_{n})_{n\ge 0}$ forms a temporally homogeneous Markov chain such that the conditional law of $(M_{n},X_{n})$ given the history up to time $n-1$ is a function of $M_{n-1}$ only (instead of $(M_{n-1},X_{n-1})$) for each $n\ge 1$. This can in fact be shown to be an equivalent property if $\fS$ is countably generated. We make the additional assumption that $(M_{n})_{n\ge 0}$ has a stationary distribution $\pi$ and  is ergodic under $\Prob_{\pi}:=\int_{\cS}\Prob(\cdot|M_{0}=s)\,\pi(\dd s)$. By ergodic decomposition, the latter assumption does not entail any loss of generality regarding our results. Due to the special Markovian structure, this entails the same for $(M_{n},X_{n})_{n\ge 0}$. Let $\mu$ be the pertinent stationary law and put $\Prob_{\mu}:=\int_{\cS\times\R^{d}}\Prob(\cdot|M_{0}=s,X_{0}=x)\,\mu(\dd  s\times\dd x)$. Notice that the laws of $(M_{n},X_{n})_{n\ge 1}$ under $\Prob_{\mu}$ and $\Prob_{\pi}$ are identical.

\vspace{.1cm}
Defining $S_{0}:=0$ and $S_{n}:=\sum_{i=1}^{n}X_{i}$ for $n\ge 1$, the bivariate sequence $(M_{n},S_{n})_{n\ge 0}$ and also $(S_{n})_{n\ge 0}$ alone are called \emph{Markov random walk (MRW)} and $(M_{n})_{n\ge 0}$ its \emph{driving} or \emph{modulating chain}. Whenever these objects are studied in stationary regime, that is, under $\Prob_{\mu}$, we may also consider a doubly infinite stationary extension $(M_{n},X_{n})_{n\in\Z}$ with associated doubly infinite random walk
\begin{align}\label{eq:doubly infinite MRW}
S_{n}\ =\ \begin{cases} \sum_{i=1}^{n}X_{i}&\text{if }n\ge 1,\\ 0,&\text{if }n=0,\\ -\sum_{i=n+1}^{0}X_{i}&\text{if }n<0. \end{cases}
\end{align}
In this framework, we call both $(M_{n},S_{n})_{n\in\Z}$ and $(M_{n},X_{n})_{n\in\Z}$
\begin{itemize}\itemsep2pt
\item \emph{null-homolo\-gous} if there exists a stationary sequence $(Y_{n})_{n\in\Z}$ such that
\begin{gather*}
X_{n}\ =\ Y_{n}-Y_{n-1}\quad\Prob_{\mu}\text{-a.s.}
\shortintertext{and thus}
S_{n}\ =\ \sign(n)(Y_{n}-Y_{0})\quad\Prob_{\mu}\text{-a.s.}
\end{gather*}
for all $n\in\Z$;\footnote{Thus, $(X_{n})_{n\in\Z}$ defined in \eqref{def-nullhom-measure1} forms a stationary sequence with associated MRW $(M_{n},S_{n})_{n\in\Z}$ under any $T$-invariant probability measure on $\Omega$, where $M_{n}:=T^{n}$, and null-homology of $\Prob$ in Definition \ref{def-nullhom-measure2} is equivalent to the null-homology of the MRW $(M_{n},S_{n})_{n\in\Z}$ under $\Prob$.}
\item \emph{strict-sense null-homologous} if, furthermore, there exists a measurable function $\xi:\cS\to\R^{d}$ such that $Y_{n}=\xi(M_{n})$ $\Prob_{\mu}$-a.s.~for all $n\ge 0$.
\end{itemize}
For the unilateral sequences $(M_{n},S_{n})_{n\ge 0}$ and $(M_{n},X_{n})_{n\ge 0}$, the stated conditions must naturally hold for all $n\ge 1$. On the other hand, these conditions persist under doubly infinite stationary extension whence it actually makes no difference whether the unilateral or the doubly infinite framework is adopted. Regarding the relation between null-homology and the formally stronger strict-sense null-homology, it is not obvious and therefore shown as part of our main result, Theorem \ref{thm-main 1} below, that they are actually equivalent.

\vspace{.1cm}
We further note that the function $\xi$ must be unique up to translation by some $c\in\R^{d}$. Namely, if $X_{n}=\zeta(M_{n})-\zeta(M_{n-1})$ $\Prob_{\mu}$-a.s.~for all $n\ge 1$ and another function $\zeta:\cS\to\R^{d}$, then 
$$ \xi(M_{0})-\zeta(M_{0})\ =\ \xi(M_{0})-\zeta(M_{0})\ =\ \ldots\quad\Prob_{\mu}\text{-a.s.} $$
In other words, $\xi-\zeta$, viewed as a function on $\cS^{\N_{0}}$, is a.s.~shift-invariant with respect to the law of $(M_{n})_{n\ge 0}$ under $\Prob_{\mu}$ and therefore as claimed equal to some $c\in\R^{d}$ by ergodicity.

\vspace{.1cm}
Null-homology arises quite naturally in connection with the lattice-type of one-dimensional MRW's. Let $(M_{n})_{n\ge 0}$ be ergodic with unique stationary law $\pi$. Following Shurenkov \cite{Shurenkov:84}, the MRW $(M_{n},S_{n})_{n\ge 0}$ is called $d$-arithmetic if $d$ is the maximal positive number such that
$$ \Prob_{\pi}\big(X_{1}\in \xi(M_{1})-\xi(M_{0})+d\Z\big)\ =\ 1 $$
for a suitable function $\xi:\cS\to [0,d)$, called \emph{shift function}. If no such $d$ exists, it is called nonarithmetic. Equivalently, $(M_{n},S_{n})_{n\ge 0}$ is $d$-arithmetic if $d>0$ is the maximal number such that $(M_{n},X_{n}-X_{n}')_{n\in\Z}$ is Markov-modulated and strict-sense null-homologous for a sequence of $d\Z$-valued random variables $(X_{n}')_{n\in\Z}$.
Namely, with $\xi$ denoting the shift function,
$$ X_{n}'\,:=\,X_{n}-\xi(M_{n})+\xi(M_{n-1}) $$ 
for $n\in\Z$.

\vspace{.1cm}
Since, given an arbitrary stationary sequence $(X_{n})_{n\in\Z}$, any of
$$ M_{n}\ :=\ (X_{n+i})_{i\in\Z},\ (X_{n-i})_{i\le 0},\ \text{or }(X_{n+i})_{i\ge 0},\quad n\in\Z $$
constitutes a modulating stationary Markov chain on a Polish state space (ergodic iff $(X_{n})_{n\in\Z}$ is ergodic and also Fellerian), we see that null-homology for stationary processes can also be studied within the framework of Markov-modulation, with some freedom to choose the driving chain. Moreover, the conditional law of each $X_{n}$ given the driving chain is always degenerate under the above choices. 

\subsection{Main results relating strict-sense null-homology and tightness} 

By definition, null-homology implies stationarity of the sequence $(S_{n}+Y_{0})_{n\in\Z}$ and thus ``almost stationarity'' as well as tightness of the random walk $(S_{n})_{n\in\Z}$ itself. Regarding tightness, the converse is established by the next result together with the nonobvious fact that null-homology and strict-sense null-homology are in fact the same.

\begin{Theorem}\label{thm-main 1}
Let $(M_{n},X_{n})_{n\in\Z}$ be a doubly infinite stationary Markov-modulated sequence of $\cS\times\R^{d}$-valued random variables with ergodic driving chain. Then the following assertions are equivalent for the associated MRW $(M_{n},S_{n})_{n\in\Z}$ defined by \eqref{eq:doubly infinite MRW}:
\begin{itemize}\itemsep2pt
\item[(a)] $(M_{n},S_{n})_{n\in\Z}$ is strict-sense null-homologous.
\item[(b)] $(M_{n},S_{n})_{n\in\Z}$ is null-homologous and $X_{n}=g(M_{n-1},M_{n})$ a.s.~for each $n\in\Z$ and some measurable $g:\cS\times\cS\to\R^{d}$.
\item[(c)] $(M_{n},S_{n})_{n\in\Z}$ is null-homologous.
\item[(d)] $(S_{n})_{n\in\Z}$ is tight.
\item[(e)] $(S_{n})_{n\ge 0}$ is tight.
\end{itemize}
\end{Theorem}

As strict-sense null-homology by definition ensures $X_{n}=\xi(M_{n})-\xi(M_{n-1})$ a.s.~for some $\xi$, we see that (a) trivially implies (b) and particularly that the $X_{n}$ given their driving chain are all a.s.~constant. In other words, the modulation is rigid in the sense of no extra randomness beyond the contribution of the driving chain. Our result may therefore be summarized by saying that this rigidity in combination with null-homology is necessary and sufficient for the tightness and thus atypical behavior of a MRW.

\subsection{Null-recurrent driving chain}

For a MRW $(M_{n},S_{n})_{n\ge 0}$ with null-re\-current driving chain, the notion of null-homology can be generalized as follows: Let $\pi$ denote the essentially unique $\sigma$-finite stationary measure of $(M_{n})_{n\ge 0}$ and $\fS_{\pi}=\{A\in\fS:0<\pi(A)<\infty\}$. Putting $\pi_{A}:=\pi(\cdot\cap A)/\pi(A)$, it is well-known that the driving chain is \emph{cycle-stationary under} $\Prob_{\pi_{A}}$ for any $A\in\fS_{\pi}$ in the sense that the sequence of cycles
$$ (M_{\tau_{n}(A)},\ldots,M_{\tau_{n+1}(A)-1}),\quad n\ge 0 $$
is stationary under $\Prob_{\pi_{A}}$ in the ordinary sense. Here the $\tau_{n}(A)$ for $n\ge 0$ denote the successive return times to $A$ of the chain with $\Prob_{\pi_{A}}(\tau_{0}(A)=0)=1$. By Markov modulation, the sequence $(M_{n},X_{n})_{n\in\Z}$ is then cycle stationary under $\Prob_{\mu_{A}}$, where $\mu_{A}$ equals the law of $(M_{\tau_{1}(A)},X_{\tau_{1}(A)})$ under $\Prob_{\mu_{A}}$. Null-homology and strict-sense null-homology for $(M_{n},X_{n})_{n\ge 0}$ and its associated MRW are now defined exactly as before, but with $(Y_{n})_{n\ge 0}$ being cycle-stationary under $\Prob_{\mu_{A}}$ for any $A\in\fS_{\pi}$. The following result provides the counterpart of Theorem \ref{thm-main 1} in the null-recurrent situation.

\begin{Theorem}\label{thm-main 2}
Let $(M_{n},X_{n})_{n\ge 0}$ be a Markov-modulated sequence with null-recurrent driving chain and associated MRW $(M_{n},S_{n})_{n\ge 0}$. Then the following assertions are equivalent:
\begin{itemize}\itemsep2pt
\item[(a)] $(M_{n},S_{n})_{n\ge 0}$ is strict-sense null-homologous.
\item[(b)] $(M_{n},S_{n})_{n\ge 0}$ is null-homologous and $X_{n}=g(M_{n-1},M_{n})$ $\Prob_{\mu_{A}}$-a.s.~for each $n\in\N$, $A\in\fS_{\pi}$ and some measurable $g:\cS\times\cS\to\R^{d}$.
\item[(c)] $(M_{n},S_{n})_{n\ge 0}$ is null-homologous.
\item[(d)] $(S_{\tau_{n}(A)})_{n\ge 0}$ is tight under $\Prob_{\mu_{A}}$ for any $A\in\fS_{\pi}$.
\end{itemize}
Provided they hold and $S_{n}=\xi(M_{n})-\xi(M_{0})$ $\Prob_{\mu}$-a.s.~for all $n$, the random walk $(S_{n})_{n\ge 0}$ itself is tight under $\Prob_{\mu_{A}}$ for all $A\in\fS_{\pi}$ iff 
the pre-image of the complement of some $[-a,a]^{d}$, $a\ge 0$, under $\xi$ has finite mass with respect to $\pi$, the stationary measure of the driving chain.
\end{Theorem}

\vspace{.1cm}
We have organized the rest of the article as follows. The proof of Theorem \ref{thm-main 1} together with further results such as an $L^{p}$-version of the coboundary theorem and another characterization of strict-sense null-homology are the content of Section \ref{sec-proofs}. A collection of instances where null-homology arises as a relevant issue are discussed in Section \ref{sec-appli}- \ref{sec-homogenization}. In particular, we discuss in Section \ref{sec-polaron} another characterization of null-homology (cf. Theorem \ref{thm-ergodic}) and show
in Section \ref{sec-homogenization} that the so-called corrector, which comes into play when aiming at an almost sure central limit theorem (CLT), is strict-sense null-homologous. We also explain that this does not contradict a seemingly different statement by Gloria \cite{Gloria:14} about the same issue. 

\section{Further results and proofs of Theorem \ref{thm-main 1} and Theorem \ref{thm-main 2}}\label{sec-proofs}

Adopting the notation from the previous section and additionally assuming $(\cS,\fS)$ to be Polish, the following theorem, stated for one-sided sequences, may be viewed as a more elaborate formulation of the equivalence between (c) and (e) in Theorem \ref{thm-main 1}. Its proof, given after one more theorem, is probabilistic and relatively short. Since assertions about the additive part of a given MRW $(M_{n},S_{n})_{n\ge 0}$ on $\cS\times\R^{d}$ are not affected by the choice of the driving chain, the Polish state space assumption does not constitute a real restriction because we can always choose as driving chain $M_{n}=(X_{n+k})_{k\ge 0}$ for $n\ge 0$ which has Polish state space $((\R^{d})^{\times\N_{0}},\cB(\R^{d})^{\times\N_{0}})$ and also the Feller property.

\begin{Theorem}\label{thm-main 3}
Let $(M_{n},X_{n})_{n\ge 0}$ be a Markov-modulated sequence in a standard model with state space $\cS\times\R^{d}$, ergodic stationary distribution $\mu$ and associated zero-delayed MRW $(M_{n},S_{n})_{n\ge 0}$, thus $S_{0}=0$ and $S_{n}=X_{1}+\ldots+X_{n}$ for $n\ge 1$. Let further $Y_{n}=X_{0}+S_{n}$ for $n\ge 0$. Provided that $(\cS,\fS)$ is Polish, the following assertions are equivalent:
\begin{itemize}\itemsep2pt
\item[(a)] $(S_{n})_{n\ge 0}$ is tight under $\Prob_{\mu}$.
\item[(b)] There exists a law $\nu$ on $\cS\times\R^{d}$ such that the sequences $(M_{n},X_{n})_{n\ge 1}$ and $(M_{n},X_{n+1},Y_{n})_{n\ge 0}$ are both stationary under $\Prob_{\nu}$. Moreover,
$$ \Prob_{\nu}((M_{n},Y_{n})\in\cdot)\,=\,\nu\quad\text{and}\quad\Prob_{\nu}((M_{n+1},X_{n+1})\in\cdot)\,=\,\mu $$
for all $n\ge 0$.
\item[(c)] $(M_{n},S_{n})_{n\ge 0}$ is null-homologous.
\end{itemize}
If these assertions hold, the $Y_{n}$ can be chosen such that
\begin{equation}\label{eq:Y_n<->X_n}
Y_{n}\ =\ f(X_{n+1},X_{n+2},\ldots)\quad\text{a.s.}
\end{equation}
for all $n\in\Z$ and some measurable function $f:\R^{\times\Z}\to\R^{d}$.
\end{Theorem}

In part (b) the stationarity of $(M_{n},X_{n})_{n\ge 1}$ is in fact a consequence of the stationarity of $(M_{n},X_{n+1},Y_{n})_{n\ge 0}$. Moreover, it should be observed that $\nu$ constitutes a stationary law for the Markov chain $(M_{n},S_{n})_{n\ge 0}$.

\vspace{.1cm}
Our last theorem is concerned with another characterization of strict-sense null-homo\-logy that arises in connection with the so-called corrector in stochastic homogenization (see Section \ref{sec-homogenization}). Its proof is postponed to the end of this section. Obviously, if there exists a measurable function $\xi:\cS\to\R^{d}$ such that
\begin{equation}\label{eq:strict-sense NH}
S_{n}\ =\ \xi(M_{n})-\xi(M_{0})\quad\Prob_{s}\text{-a.s.~for all }s\in\cS\text{ and }n\ge 0,
\end{equation}
then the \emph{closed-loop condition}
\begin{equation}\label{closed-loop}
\Prob(S_{n}=0|M_{0}=M_{n})\,=\,1\quad\Prob_{s}\text{-a.s.~for all }s\in\cS
\end{equation}
holds which means that, whenever the driving chain returns to its initial state and thus completes a loop, then so does the modulated random walk. For the case when $(M_{n})_{n\ge 0}$ is an irreducible, but not necessarily recurrent Markov chain on a discrete state space, the subsequent theorem asserts that the converse is also true.

\begin{Theorem}\label{thm-main 4}
Let $(M_{n},S_{n})_{n\ge 0}$ be a MRW such that its driving chain $(M_{n})_{n\ge 0}$ is irreducible with discrete state space $\cS$. Then the following assertions are equivalent:
\begin{itemize}\itemsep2pt
\item[(a)] $(M_{n},S_{n})_{n\ge 0}$ satisfies the closed-loop condition
$$ \Prob(S_{n}=0|M_{0}=M_{n})\,=\,1\quad\Prob_{s}\text{-a.s.~all }s\in\cS. $$
\item[(b)] $(M_{n},S_{n})_{n\ge 0}$ is strict-sense null-homologous in the sense that $g(s,t)=\xi(t)-\xi(s)$ for all $s,t\in\cS$ and some $\xi:\cS\to\R^{d}$.
\end{itemize}
\end{Theorem}

\begin{proof}[Proof of Theorem \ref{thm-main 3}]
``(a)$\RA$(b)'' Defining $\wh{M}_{n}:=(M_{n+k})_{k\ge 0}$ and $\wh{X}_{n},\wh{Y}_{n}$ accordingly,
consider the Feller chain $(\wh{M}_{n},\wh{X}_{n},\wh{Y}_{n})_{n\ge 0}$ with state space $\X:=\cS^{\times\N_{0}}\times(\R^{d})^{\times\N_{0}}\times(\R^{d})^{\times\N_{0}}$ and generated by the continuous map $\Psi:\X\to\X$,
$$ ((s_{k})_{k\ge 0},(x_{k})_{k\ge 0},(y_{k})_{k\ge 0})\mapsto ((s_{k})_{k\ge 1},(x_{k})_{k\ge 1},(y_{k-1}+x_{k})_{k\ge 1}) $$
(this is where the additional assumption on $(\cS,\fS)$ enters).
Under $\Prob_{\mu}$, the chain is stationary in the first two coordinates and therefore tight as a whole if (a) holds. By Prokhorov's theorem, the relative compactness of the distributional time averages
\begin{align*}
Q_{n}\ :=\ \frac{1}{n}\sum_{j=1}^{n}\Prob_{\mu}((\wh{M}_{j},\wh{X}_{j},\wh{Y}_{j})\in\cdot),\quad n\ge 1
\end{align*}
follows and thus the existence of $1\le n_{1}<n_{2}<\ldots$ such that, as $k\to\infty$, both $Q_{n_{k}}$ and $Q_{n_{k}+1}$ converge weakly to a probability law on $\X$ which constitutes a stationary distribution of $(\wh{M}_{n},\wh{X}_{n},\wh{Y}_{n})_{n\ge 0}$ (here the Feller property enters). In particular, the trivariate Markov chain $(M_{n},X_{n},Y_{n})_{n\ge 0}$ has a stationary law $\Lambda$ satisfying $\Lambda(\cdot\times\cdot\times\R^{d})=\mu$. Moreover, its transition kernel, say $K((s,x,y),\cdot)$, does not depend on $x$ because $Y_{n}=Y_{n-1}+X_{n}$ for each $n\ge 1$ and by the Markov-modulated structure of $(M_{n},X_{n})_{n\ge 0}$. Namely, with $\Prob_{s,x,y}:=\Prob(\cdot|M_{0}=s,X_{0}=x,Y_{0}=y)$, we infer
\begin{align*}
K((s,x,y),A\times B\times C)\ &=\ \Prob_{s,x,y}(M_{1}\in A,X_{1}\in B,Y_{0}+X_{1}\in C)\\
&=\ \Prob_{s,x,y}(M_{1}\in A,X_{1}\in B,y+X_{1}\in C)\\
&=\ \Prob_{s}(M_{1}\in A,X_{1}\in B,y+X_{1}\in C)
\end{align*}
for all $A\in\fS$ and $B,C\in\cB(\R^{d})$. After these observations define
$$ \nu\,:=\,\Lambda(\cdot\times\R^{d}\times\cdot) $$
and recall that $Y_{0}=X_{0}$. It follows that
\begin{align*}
\Prob_{\nu}((M_{1},X_{1})\in\cdot)\ &=\ \int_{\cS\times\R^{d}}\Prob_{s}((M_{1},X_{1})\in\cdot)\ \nu(\dd s\times\dd x)\\
&=\ \int_{\cS}\Prob_{s}((M_{1},X_{1})\in\cdot)\ \pi(\dd s)\\
&=\ \Prob_{\pi}((M_{1},X_{1})\in\cdot)\ =\ \mu.
\end{align*}
and therefore that $(M_{n},X_{n})_{n\ge 1}$ is stationary under $\Prob_{\nu}$ as claimed. Finally, as $(M_{n},X_{n+1},Y_{n})_{n\ge 0}$ is stationary under $\Prob_{\Lambda}$, the same holds true under $\Prob_{\nu}$ because, again using that $K((s,x,y),\cdot)$ does not depend on $x$,
\begin{gather*}
\Prob_{\nu}((M_{0},X_{1},Y_{0})\in A\times B\times C)\ =\ \int_{A\times C}\Prob_{s}(X_{1}\in B)\ \nu(\dd s\times\dd y)\\
\begin{split}
&=\ \int_{A\times\R^{d}\times C}\Prob_{s}(X_{1}\in B)\ \Lambda(\dd s\times\dd x\times\dd y)\\
&=\ \Prob_{\Lambda}((M_{0},X_{1},Y_{0})\in A\times B\times C)\\
&=\ \Prob_{\Lambda}((M_{n},X_{n+1},Y_{n})\in A\times B\times C)
\end{split}
\end{gather*}
for all $n\ge 1$ and $A,B,C$ as before. In particular, we find by choosing $B=\R^{d}$, that the law of $(M_{n},Y_{n})$ under $\Prob_{\nu}$ equals $\nu$ for all $n$. This completes the proof of (b).

\vspace{.1cm}
``(b)$\RA$(c)'' Since the laws of $(M_{n},S_{n})_{n\ge 0}$ under $\Prob_{\nu}$ and $\Prob_{\mu}$ are the same (as under $\Prob_{\pi}$), the null-homology follows from $Y_{n}=X_{0}+S_{n}=Y_{0}+S_{n}$ and thus $S_{n}=Y_{n}-Y_{0}$ $\Prob_{\nu}$-a.s.~for all $n\ge 1$. 

\vspace{.1cm}
Left with the proof of \eqref{eq:Y_n<->X_n}, we begin with some preliminary remarks. If $\cI$ denotes the invariant $\sigma$-field of $\wh{Y}=(Y_{n})_{n\ge 0}$, which is generated by the random variables $g(\wh{Y})$ for shift-invariant functions $g:(\R^{d})^{\times\N}\to\R^{d}$, then any $\wh{Y}'=(Y_{n}')_{n\ge 0}$ with $Y_{n}'=Y_{n}+Z$ for an $\cI$-measurable $\R^{d}$-valued $Z$ can be used instead of $(Y_{n})_{n\ge 0}$ in the sense that $(M_{n},X_{n+1},Y_{n}')_{n\ge 0}$ is still stationary under $\Prob=\Prob_{\nu}$ and
\begin{equation}\label{eq:homology eq}
X_{n}\,=\,Y_{n}'-Y_{n-1}'\quad\Prob\text{-a.s.~for all }n\ge 1.
\end{equation}
This observation allows us to pick the $Y_{n}'$ in such a way that their conditional median given $\cI$ is 0. Namely, let $F(\omega,\cdot):=\Prob(Y_{0}\le\cdot|\cI)(\omega)$ be a regular version of the conditional distribution function of $Y_{0}$ and thus of each $Y_{n}$ given $\cI$. Their conditional median $R_{0}$ is defined as
$$ R_{0}(\omega)\ :=\ \sup\left\{r\in\R:F(\omega,r)<\frac{1}{2}\right\}. $$
and $\cI$-measurable. Consequently, a sequence with the desired property is obtained by replacing $Y_{n}$ with $Y_{n}-R_{0}$ for $n\ge 0$.

\vspace{,1cm}
Proceeding with the proof of \eqref{eq:Y_n<->X_n}, we may now assume $R_{0}=0$ a.s.~and  also $d=1$, for instance by considering the $d$ components of $X_{n}$ separately. The following argument is similar to the one given by Bradley in \cite{Bradley:95} for the nonstationary situation, but it is simpler and shorter because the sequence $(Y_{n})_{n\ge 0}$ he needs to define is here already given and stationary. In particular, we do not require Koml\'os' law of large numbers.

\vspace{.1cm}
Put $S_{k,n}:=S_{k+n}-S_{k}$ and note that $Y_{k+n}=Y_{k}+S_{k,n}$ for all $k,n\in\N_{0}$. Birkhoff's ergodic theorem ensures that, as $n\to\infty$,
$$ F_{k,n}(\omega,r)\ :=\ \frac{1}{n}\sum_{j=1}^{n}\1_{(-\infty,r]}(Y_{k+n}(\omega)) $$
converges to $F(\omega,r)$ for all $(k,r)\in\N_{0}\times\R$ and $\Prob_{\nu}$-almost all $\omega\in\Omega$, say $\omega\in\Omega_{r}$. By a standard selection argument, we can also determine a set $\Omega'$ of probability one such that $\Omega'\subset\Omega_{r}$ for all rational $r$ and
$$ F(\omega,r-)\ \le\ \liminf_{n\to\infty}F_{k,n}(\omega,r)\ \le\ \limsup_{n\to\infty}F_{k,n}(\omega,r)\ \le\ F(\omega,r) $$
for all $(\omega,r)\in\Omega'\times\R\backslash\Q$. Moreover $R_{0}=0$ on $\Omega'$. Next observe that
$$ G_{k,n}(r)\ :=\ \frac{1}{n}\sum_{j=1}^{n}\1_{(-\infty,r]}(S_{k,n}) $$
satisfies
$$ G_{k,n}(\omega,r)\ =\ F_{k,n}(\omega,r+Y_{k}(\omega)) $$
for all $k,n\in\N_{0}$. It follows for all $\omega\in\Omega'$ and $k\in\N_{0}$ that
\begin{gather*}
\lim_{n\to\infty}G_{k,n}(\omega,r)\ =\ F(\omega,r+Y_{k}(\omega))
\intertext{if $r+Y_{k}(\omega)$ is rational, and}
\begin{split}
F(\omega,r+Y_{k}(\omega)-)\ &\le\ \liminf_{n\to\infty}G_{k,n}(\omega,r)\\
&\le\ \limsup_{n\to\infty}G_{k,n}(\omega,r)\ \le\ F(\omega,r+Y_{k}(\omega))
\end{split}
\end{gather*}
for all $r\in\R$. Finally defining the stationary sequence
\begin{align*}
Z_{k}\ :=\ \sup\left\{r\in\Q:\limsup_{n\to\infty}G_{k,n}(r)<\frac{1}{2}\right\}
\end{align*}
for $k\in\N_{0}$, we have $Z_{k}=f(X_{k+1},X_{k+2},\ldots)$ for some measurable function $f:\R^{\times\N}\to\R\cup\{\pm\infty\}$. Furthermore $Z_{k}=-Y_{k}$ a.s.~because, if $\omega\in\Omega'$, then
\begin{gather*}
\limsup_{n\to\infty}G_{k,n}(\omega,r-Y_{k}(\omega))\ =\ \limsup_{n\to\infty}F_{k,n}(\omega,r)\ \ge\ F(r-)\ \ge\ \frac{1}{2}
\shortintertext{for $r>0\equiv R_{0}$, and}
\limsup_{n\to\infty}G_{k,n}(\omega,r-Y_{k}(\omega))\ =\ \limsup_{n\to\infty}F_{k,n}(\omega,r)\ \le\ F(r)\ <\ \frac{1}{2}
\end{gather*}
for $r<0$. Hence $Y_{k}=f(X_{k+1},X_{k+2},\ldots)$ a.s.~and $f$ is real-valued.
\end{proof}

Before proceeding to the remaining proofs, we insert the $L^{p}$-version of the coboundary theorem, with a short proof and for all $p>0$. It goes back to Leonov \cite{Leonov:61} for $p=2$ and even wide-sense stationary sequences, to Lalley \cite{Lalley:86} for $p=1$ (as already mentioned), and it is stated by Aaronson and Weiss \cite{AarWeiss:00} for general $p\ge 1$.

\begin{Cor}\label{cor-main 3}
In the situation of the previous theorem, let $|\cdot|$ denote Euclidean norm on $\R^{d}$. Then the following assertions are equivalent for any $p>0$:
\begin{itemize}\itemsep2pt
\item[(a)] $\sup_{n\ge 0}\Erw_{\mu}|S_{n}|^{p}<\infty$.
\item[(b)] $\Erw_{\nu}|Y_{0}|^{p}=\int|x|^{p}\,\nu(dx)<\infty$.
\end{itemize}
We call $(S_{n})_{n\ge 0}$ $L^{p}$-null-homologous under these conditions.
\end{Cor}

\begin{proof}
Plainly, we must only verify that (a) implies (b). Returning to the proof of Theorem \ref{thm-main 3}, note that $\sup_{n}\Erw_{\mu}|S_{n}|^{p}<\infty$ implies
$$ \sup_{n\ge 1}\int|x|^{p}\ Q_{n}(dx)\ <\ \infty. $$
Choosing the $n_{k}$ such that $Q_{n_{k}}\weakly\nu$, we can always define, on a possibly enlarged probability space, a sequence $(W_{k})_{k\ge 0}$ such that $\Prob_{\mu}(W_{0}\in\cdot)=\nu$, $\Prob_{\mu}(W_{k}\in\cdot)=Q_{n_{k}}$ for $k\ge 1$ and $W_{k}\to W_{0}$ $\Prob_{\mu}$-a.s. Now
$$ \Erw_{\nu}|Y_{0}|^{p}\ =\ \Erw_{\mu}|W_{0}|^{p}\ \le\ \liminf_{k\to\infty}\Erw_{\mu}|W_{k}|^{p}<\infty $$
by Fatou's lemma.
\end{proof}

\begin{proof}[Proof of Theorem \ref{thm-main 1}]
Since ``(a)$\RA\ldots\RA$(e)'' is straightforward and ``(e)$\RA$(c)'' is ensured by Theorem \ref{thm-main 3}, only ``(c)$\RA$(a)'' remains to be verified. By the afore-mentioned result, we have that $S_{n}=Y_{n}-Y_{0}$ a.s.~for all $n\in\Z$ and a stationary sequence $(Y_{n})_{n\in\Z}$ satisfying $Y_{n}=f(X_{n+1},X_{n+2},\ldots)$ a.s.~for all $n$ and some measurable $f:(\R^{d})^{\times\N}\to\R^{d}$. We recall the assumption that the driving chain $(M_{n})_{n\in\Z}$ is ergodic and also exclude the trivial case when all $X_{n}$ and thus all $Y_{n}$ vanish almost surely.

\vspace{.1cm}
As a first step, we show that the conditional law of $S_{n}$ given $M_{0},M_{n}$ must be a.s.~degenerate for each $n\in\N$, thus $S_{n}=g_{n}(M_{0},M_{n})$ a.s.~for some measurable $g_{n}:\cS\times\cS\to\R^{d}$.

\vspace{.1cm}
For $x,y\in\cS$, $n\in\Z$ and $t\in\R^{d}$, let 
$$ \vph^{x}(t)=\Erw(e^{itY_{0}}|M_{0}=x)\quad\text{and}\quad\psi_{n}^{x,y}(t)=\Erw(e^{itS_{n}}|M_{0}=x,M_{n}=y) $$ 
be the conditional Fourier transforms (FT) of $Y_{0}$ given $M_{0}=x$ and of $S_{n}$ given $M_{0}=x,\,M_{n}=y$, respectively. Using the Markov-modulated structure, we infer that the conditional FT of $S_{kn}$ given $M_{0},M_{kn}$ satisfies
$$ \psi_{kn}^{M_{0},M_{kn}}(t)\ =\ \Erw\left[\prod_{j=1}^{k}\psi_{n}^{M_{(j-1)n},M_{jn}}(t)\bigg|M_{0},M_{kn}\right]\quad\text{a.s.} $$
for all $k,n\in\N$. Next, if there exist $m\in\N$ and a subset $C$ of $\cS^{2}$ such that
$\Prob((M_{0},M_{m})\in C)>0$ and the conditional law of $S_{n}$ given $M_{0}=x,M_{m}=y$ is nondegenerate for $(x,y)\in C$, then $C$ may in fact be chosen in such a way that, for some $t_{0}>0$,
$$ |\psi_{m}^{M_{0},M_{m}}(t_{0})|\,<\,1\quad\text{a.s.~on }\{(M_{0},M_{m})\in C\} $$
holds together with $\Erw|\vph^{M_{0}}(t_{0})|>0$. Birkhoff's ergodic theorem then further implies that
$$ \frac{1}{k}\sum_{j=1}^{k}\log|\psi_{m}^{M_{(j-1)m},M_{jm}}(t_{0})|\ \xrightarrow{k\to\infty}\ \Erw\log|\psi_{m}^{M_{0},M_{m}}(t_{0})|\ <\ 0\quad\text{a.s.} $$
Consequently,
\begin{align*}
|\psi_{km}^{M_{0},M_{km}}(t_{0})|\ &=\ \left|\Erw\left[\prod_{j=1}^{k}\psi_{m}^{M_{(j-1)m},M_{jm}}(t_{0})\Bigg|M_{0},M_{km}\right]\right|\\
&\le\ \Erw\left[\exp\left(\sum_{j=1}^{k}\log|\psi_{m}^{M_{(j-1)m},M_{jm}}(t_{0})|\right)\Bigg|M_{0},M_{km}\right]\\
&=\ \Erw\left[\exp\Big(k\,\Erw\log|\psi_{m}^{M_{0},M_{m}}(t_{0})|(1+o(1))\Big)\Bigg|M_{0},M_{km}\right]\\
&\xrightarrow{k\to\infty}\ 0\quad\text{a.s.}
\end{align*}

On the other hand, use $Y_{0}=Y_{n}-S_{n}$ and the conditional independence of $S_{n}$ and $Y_{n}=f(X_{n+1},X_{n+2},\ldots)$ given $M_{0},M_{n}$ to obtain the equation
\begin{align*}
\vph^{M_{0}}(t)\ =\ \Erw\Big[\vph^{M_{n}}(t)\ovl{\psi_{n}^{M_{0},M_{n}}(t)}\Big|M_{0}\Big]\quad\text{a.s.}
\end{align*}
for any $n\in\N$ (with $\ovl{z}$ denoting the complex conjugate of $z$). But for $t=t_{0}$, this provides with the help of the dominated convergence theorem that
\begin{align*}
0\ <\ \Erw|\vph^{M_{0}}(t_{0})|\ &=\ \lim_{k\to\infty}\Erw\Big|\vph^{M_{km}}(t_{0})\ovl{\psi_{km}^{M_{0},M_{km}}(t_{0})}\Big|\\
&\le\ \Erw\left[\lim_{k\to\infty}\Big|\psi_{km}^{M_{0},M_{km}}(t_{0})\Big|\right]\ =\ 0
\end{align*}
which is impossible.

\vspace{.1cm}
Having verified that $S_{n}=g_{n}(M_{0},M_{n})$ a.s.~for all $n\in\N$ and a suitable function $g_{n}$, we now see that
$$ Y_{n}\ =\ h(M_{n},M_{n+1},\ldots)\quad\text{a.s.} $$
for all $n$ and some measurable $h:\cS^{\times\N}\to\R^{d}$. It follows that
$$ h(M_{0},M_{1},\ldots)\ =\ Y_{0}\ =\ Y_{n}-S_{n}\ =\ H_{n}(M_{0},M_{n},M_{n+1},\ldots)\quad\text{a.s.} $$
where $H_{n}(x,y_{0},y_{1},\ldots):=f(y_{1},y_{2},\ldots)-g_{n}(x,y_{0})$, that is, $h(M_{0},M_{1},\ldots)$ is a.s.~constant in $M_{1},\ldots,M_{n-1}$. But this being true for each $n\ge 2$ means that $h(M_{0},M_{1},\ldots)$ must be a.s.~constant in $M_{1},M_{2},\ldots$ and thus a.s.~equal to $\xi(M_{0})$ for some $\xi:\cS\to\R^{d}$. This completes the proof.
\end{proof}

\begin{proof}[Proof of Theorem \ref{thm-main 2}]
We must only prove ``(d)$\RA$(a)'' as one can readily see. Put $(M_{n}^{A},X_{n}^{A}):=(M_{\tau_{n}(A)},S_{\tau_{n}(A)}-S_{\tau_{n-1}(A)})$ for $n\ge 1$ and $(M_{0}^{A},X_{0}^{A}):=(M_{0},X_{0})$. For $A\in\fS_{\pi}$, this sequence is Markov-modulated, ergodic with unique stationary law $\mu_{A}$ and thus stationary under $\Prob_{\mu_{A}}$. The assumed tightness of the associated MRW $(S_{\tau_{n}(A)})_{n\ge 0}$ implies by Theorem \ref{thm-main 1} that it is strict-sense null-homolo\-gous. So $X_{n}^{A}=\xi(M_{n}^{A})-\xi(M_{n-1}^{A})$ $\Prob_{\mu_{A}}$-a.s.~for all $n\ge 1$ and a suitable function $\xi_{A}:A\to\R^{d}$ which is unique up to a translation by some $c\in\R^{d}$ as noted after the definition of strict-sense null-homology.

\vspace{.1cm}
Since $\pi$ is $\sigma$-finite, we can pick a strictly increasing sequence $(\cS_{m})_{m\ge 1}$ such that $\cS_{m}\in\fS_{\pi}$ for each $m$ and $\bigcup_{m}\cS_{m}=\cS$. As just pointed out,  we have $S_{n}^{\cS_{m}}=\xi_{m}(M_{n}^{\cS_{m}})-\xi_{m}(M_{0}^{\cS_{m}})$ $\Prob_{\mu_{\cS_{m}}}$-a.s. for all $m,n\ge 1$ and functions $\xi_{m}:\cS_{m}\to\R^{d}$, and we show now that these functions can be chosen such that the restriction of $\xi_{m+1}$ to $\cS_{m}$ equals $\xi_{m}$ for each $m$. First note that, by the nested structure of the $\cS_{m}$, the $S_{n}^{\cS_{m}}$, $n\ge 1$, form a subsequence of $(S_{n}^{\cS_{m+1}})_{n\ge 1}$. Consequently,
\begin{gather*}
S_{n}^{\cS_{m}}\ =\ \xi_{m}(M_{\tau_{n}(\cS_{m})})-\xi_{m}(M_{0})\ =\ \xi_{m+1}(M_{\tau_{n}(\cS_{m})})-\xi_{m+1}(M_{0})
\shortintertext{and thus}
\xi_{m+1}(M_{0})-\xi_{m}(M_{0})\ =\ \xi_{m+1}(M_{n}^{\cS_{m}})-\xi_{m}(M_{n}^{\cS_{m}})
\end{gather*}
holds $\Prob_{\mu_{\cS_{m}}}\text{-a.s.}$ for all $m,n\in\N$. This further yields
$$ \xi_{m+1}(M_{0})-\xi_{m}(M_{0})\ =\ \lim_{n\to\infty}\frac{1}{n}\sum_{k=1}^{n}\Big(\xi_{m+1}(M_{n}^{\cS_{m}})-\xi_{m}(M_{n}^{\cS_{m}})\Big) $$
$\Prob_{\mu_{\cS_{m}}}\text{-a.s.}$ for any $m$. But the right-hand limit is measurable with respect to the invariant $\sigma$-field of the ergodic and $\Prob_{\mu_{\cS_{m}}}$-stationary Markov chain $(M_{n}^{\cS_{m}})_{n\ge 0}$, hence $\Prob_{\mu_{\cS_{m}}}\text{-a.s.}$ equal to some $c_{m}\in\R^{d}$. But the $\xi_{m}$ being unique only up to translation, they may be chosen such that, for all $m\ge 1$, $\xi_{m+1}(s)=\xi_{m}(s)$ a.s.~with respect to
$\Prob_{\mu_{\cS_{m}}}(M_{0}\in\cdot)=\pi_{\cS_{m}}$. In other words, there exists a mapping $\xi:\cS\to\R^{d}$ that $\Prob_{\mu}$-a.s.~equals $\xi_{m}$ on $\cS_{m}$ for any $m$ and
$$ S_{\tau_{1}(\cS_{m})}\ =\ S_{1}^{\cS_{m}}\ =\ \xi(M_{\tau_{1}(\cS_{m})})-\xi(M_{0})\quad\Prob_{\mu}\text{-a.s.} $$
But $\tau_{1}(\cS_{m})\to 1$ $\Prob_{\mu}$-a.s.~as $m\to\infty$ finally proves $X_{1}=\xi(M_{1})-\xi(M_{0})$ and thus, by stationarity, $X_{n}=\xi(M_{n})-\xi(M_{n-1})$ $\Prob_{\mu}$-a.s.~for any $n$.
\end{proof}

\begin{proof}[Proof of Theorem \ref{thm-main 4}]
We must only prove ``(a)$\RA$(b)''. Put $p_{xy}^{n}:=\Prob_{x}(M_{n}=y)$ and let $\psi_{n}^{x,y}$ be as in the previous proof for $x,y\in\cS$ and $n\ge 0$. Since $(M_{n})_{n\ge 0}$ is irreducible, $p_{xy}^{n(x,y)}>0$ for all $x,y\in\cS$ and some minimal integer $n(x,y)$. The closed-loop condition implies, with $n=n(x,y)+n(y,x)$,
\begin{gather*}
\Prob(S_{n}=0|M_{0}=M_{n}=x)\,=\,1
\shortintertext{and thereupon}
1\ =\ \psi_{n}^{x,x}(t)\ =\ \Erw\Big(\psi_{n(x,y)}^{x,M_{n(x,y)}}(t)\psi_{n(y,x)}^{M_{n(x,y)},x}(t)\Big|M_{0}=M_{n}=x\Big)
\end{gather*}
which in turn yields $\psi_{n(x,y)}^{x,s}(t)\psi_{n(y,x)}^{s,x}(t)=1$ for all $s\in\cS$ (including $s=y$) such that $\Prob(M_{n(x,y)}=s|M_{0}=M_{n}=x)>0$. This shows that 
$$ \psi_{n(x,y)}^{x,y}(t)=e^{itg(x,y)}=e^{-itg(y,x)}=\ovl{\psi_{n(y,x)}^{y,x}(t)}$$
for all $t\in\R^{d}$ and some $g:\cS\times\cS\to\R^{d}$ which satisfies $g(x,y)=-g(y,x)$. In particular, $X_{n}=g(M_{n-1},M_{n})$ $\Prob_{s}\text{-a.s.}$ for all $s\in\cS$ and $n\in\N$.
Next, a similar argument using a path $x\to y\to z\to x$ of length $n=n(x,y)+n(y,z)+n(z,x)$ provides us with
$$ g(x,y)+g(y,z)+g(z,x)\ =\ g(x,y)+g(y,z)-g(x,z)\ =\ 0 $$
so that, fixing a reference state $x$ and defining $\xi(y):=g(x,y)$ for all $y\in\cS$ (thus $\xi(x)=0$), we arrive at the desired conclusion $g(y,z)=\xi(z)-\xi(y)$ for all $y,z\in\cS$.
\end{proof}

\section{Null-homology in applications}\label{sec-appli}

This section is devoted to a small collection of instances where null-homology arises as a relevant or even crucial issue.

\subsection{Fluctuation theory for Markov random walks}\label{subsec-3.1}

Fluctuation theory for Markov random walks aims to extend results about the fine structure of ordinary random walks with i.i.d.~increments (like fluctuation-type, recurrence versus transience, ladder variables, arcsine law) to Markov random walks \cite{Alsmeyer:00,Alsmeyer:01,AlsBuck:17c,AlsBuck:19}. For instance, the well-known trichotomy that an ordinary random walk $(S_{n})_{n\ge 0}$ must be either positive divergent $(S_{n}\to\infty\text{ a.s.})$, negative divergent $(S_{n}\to-\infty\text{ a.s.})$, or oscillating $(\liminf_{n\to\infty}S_{n}=-\infty$ and $\limsup_{n\to\infty}S_{n}=\infty\text{ a.s.})$ embarks on the simple but basic fact that the only exceptional case is the trivial one when the increments are a.s.~zero. In the simplest Markov-modulated case, where the driving chain is positive recurrent on a finite or countable state space, the same trichotomy only holds when ruling out the class of all strict-sense null-homologous MRWs, see \cite[Section 5, especially Prop.~5.4]{AlsBuck:17c}. In other words, the exceptional class is no longer just a singleton.

\vspace{.1cm}
Null-homology also arises in connection with the ladder epochs of a MRW $(M_{n},S_{n})_{n\ge 0}$. To be more precise, let the driving chain be positive Harris recurrent with stationary law $\pi$ and the stationary drift $\Erw_{\pi}S_{1}$ of the random walk be positive. Define $\sigma_{0}:=0$ and
$$ \sigma_{n}\,:=\,\inf\{k>\sigma_{n-1}:S_{k}>S_{\sigma_{n-1}}\} $$
for $n\ge 1$. Then $(M_{\sigma_{n}},\sigma_{n})_{n\ge 0}$ is also a MRW with the same lattice-span as the MRW $(M_{n},n)_{n\ge 0}$ that has deterministic additive part \cite[Thm.~1(ii)]{Alsmeyer:00}. Moreover, if $(M_{n})_{n\ge 0}$ is $d$-periodic for some $d>0$ with $\Prob_{\pi}$-a.s.~unique cyclic classes $\cS_{0},\ldots,\cS_{d-1}$ (indexed in correct transitional order), then $(M_{n},n)_{n\ge 0}$ is $d$-arithmetic with shift function
$$ \gamma:\cS\to\{0,\ldots,d-1\},\quad x\ \mapsto\ \sum_{r=0}^{d-1}(d-r)\1_{\cS_{r}}(x), $$
cf.~\cite[Section 4]{Alsmeyer:00}. In other words, the ladder epochs $\sigma_{n}$ are $d$-arithmetic up to the null-homologous sequence $(\gamma(M_{\sigma_{n}})-\gamma(M_{\sigma_{n-1}}))_{n\ge 1}$.

\subsection{Null-homologous $L^{2}$, but not $L^{2}$ null-homologous}

A strongly mixing stationary sequence of nondegenerate, zero-mean and pairwise uncorrelated random variables such that the associated RW fails to satisfy the central limit theorem (CLT) was given by Herrndorf \cite{Herrndorf:83} who states as part of his result that this RW is tight and $\gamma:=\min_{n\ge 1}\Prob(S_{n}=0)>0$. So it is null-homologous, and if $S_{n}=Y_{n}-Y_{0}$ a.s.~for stationary $Y_{0},Y_{1},\ldots$,  then furthermore $\Prob(Y_{n}=Y_{0})\ge\gamma$ for all $n\ge 1$. On the other hand, since $\Erw S_{n}^{2}=n\,\Erw X_{1}^{2}>0$ for all $n$ and thus $\sup_{n\ge 1}\Erw S_{n}^{2}=\infty$, the $Y_{n}$ cannot be square-integrable by Corollary \ref{cor-main 3}. This was also pointed out in \cite{Bradley:95} and shows that an $L^{2}$-sequence can be null-homologous, but not $L^{2}$ null-homologous.

\subsection{Poisson equation}\label{subsec-3.3}

Let $(M_{n})_{n\ge 0}$ be an ergodic Markov chain with transition kernel $P$ and stationary law $\pi$. Then any pair $(f,g)$ of real-valued $L^{1}(\pi)$-functions which satisfies the Poisson equation $g=f+Pg$ can be associated with a RW, namely
\begin{equation}\label{RW Poisson equation}
S_{n}(f)\ :=\ \sum_{k=1}^{n}f(M_{k}),\quad n\ge 0,
\end{equation}
which, under $\Prob_{\pi}$, is also a martingale with stationary increments up to a null-homologous sequence. To see this, just notice that $S_{n}(f)=W_{n}+R_{n}$ with martingale part
\begin{gather*}
W_{n}\,:=\,\sum_{k=1}^{n}\big(g(M_{k})-Pg(M_{k-1})\big)
\shortintertext{and null-homologous part}
R_{n}\,:=\,Pg(M_{0})-Pg(M_{n})
\end{gather*}
for $n\ge 0$. Provided that $f,g\in L^{2}(\pi)$, Gordin and Lif{\v s}ic \cite{GordinLifsic:78} showed that this allows to derive a CLT for $n^{-1/2}S_{n}(f)$. Namely, since null-homology implies $n^{-1/2}R_{n}\to 0$ in probability, the problem reduces to a CLT for the martingale part which may be found e.g.~in \cite{Helland:82}. The same approach was used by Benda \cite{Benda:98} and by Wu and Woodroofe \cite{WuWood:00} to establish CLT's for certain contractive iterated function systems. It was pointed out by Woodroofe \cite{Woodroofe:92} that, given $f\in L^{2}(\pi)$, a solution $g$ to Poisson's equation does indeed exist if
$$ g_{n}\,:=\,\sum_{k=1}^{n}P^{k}f\ \xrightarrow{n\to\infty}\ g\quad\text{in }L^{2}(\pi), $$
and that this condition is also necessary if the doubly infinite extension of $(M_{n})_{n\ge 0}$ has trivial left tail $\sigma$-field.

\vspace{.1cm}
When there is no solution $g$, a perturbed version of the Poisson equation, viz. $(1+\eps)g_{\eps}=f+g_{\eps}$ for $\eps>0$, can be considered instead, which has unique solution
$g_{\eps}=\sum_{n\ge 1} \frac{P^{n-1}f}{(1+\eps)^{n}}\in L^{2}(\pi)$ if $f\in L^{2}(\pi)$. This was done by Kipnis and Varadhan \cite{KipVara:86} (for reversible $P$, see also the next subsection) and by Maxwell and Woodroofe \cite{MaxwellWood:00}. Further extensions of the CLT for RWs $S_{n}(f)$ as defined in \eqref{RW Poisson equation} were obtained by Derriennic and Lin \cite{DerriennicLin:01,DerriennicLin:01b,DerriennicLin:03} by building upon a fractional version of the Poisson equation which actually also leads to a fractional notion of null-homology (by them called fractional coboundary). We refrain from giving further details.

\section{Null-homology in the Polaron problem}\label{sec-polaron}

Let $\Omega_0$ denote the space of continuous functions $\omega:\mathbb R\to\mathbb R^d$ vanishing at the origin.  Then for any $t\in \R$, we have a shift $\theta_t:\Omega_0\to \Omega_0$ defined by $(\theta_t\omega)(\cdot)=\omega(t+\cdot)-\omega(\cdot)$, and 
we denote by $\mathcal M_{\mathrm{si}}(\Omega_0)$  the space of $\theta_t$-invariant probability measures on $\Omega_0$, or the space of {\it processes with stationary increments}. We also have an action on $\Omega_0\otimes \mathbb R^{d}$, with a slight abuse of notation again denoted $\theta_{t}$, that is defined by $\theta_t(\omega,x)=(\omega(t+\cdot)-\omega(t), x+ \omega(t))$. Then we denote by $\mathcal M_{\mathrm{s}}(\Omega_0\otimes \R^d)$ the space of $\theta_t$-invariant probability measures on $\Omega_0\otimes\R^d$, or the space of {\it stationary processes}. 

\vspace{.1cm}
In this context, the issue of null-homology appeared in \cite{MV:18} while identifying the strong coupling limit of Polaron path measures. Note that not every process with stationary increments appear as the {\it increments} of another stationary process. In the absence of tightness, a general criterion was provided in \cite[Theorem 3.1]{MV:18} which determined when any $\mathbb Q\in \mathcal M_{\mathrm{si}}(\Omega_0)$ admits a cocycle representation -- that is, when any $\mathbb Q\in \mathcal M_{\mathrm{si}}(\Omega_0)$ appears as the increments of some $\mathbb Q^\prime\in \mathcal M_{\mathrm{s}}(\Omega_0\otimes\mathbb R^d)$): 

\begin{Theorem}\label{thm-ergodic}(\cite[Theorem 3.1]{MV:18})
Let $\beta$ be an ergodic process with\nobreak\ stat\-ionary increments, i.e., $\beta\in \mathcal M_{\mathrm{si}}(\Omega_0)$ is a $\theta_t$-invariant and ergodic probability distribution on $\Omega_0$. Then either
\begin{equation}\label{eq:cocycle}
\lim_{\eps\to 0}\mathbb E^{\beta}\bigg[\eps \int_0^\infty \mathrm e^{-\eps t} V(\omega(t)-\omega(0)) \dd t\bigg]=0
\end{equation}
for all continuous  functions $V: \mathbb R^d\to \R$ with $\lim_{|x|\to\infty} |V(x)|=0$
or there is a $\theta_t$-invariant probability distribution $\mathbb Q\in \mathcal M_{\mathrm s}(\Omega_0\otimes\R^d)$ such that ${\beta}={\mathbb Q}^{\ssup 1}$,\footnote{Here, $\mathbb Q^{\ssup 1}$ stands for the ``first marginal" of $\mathbb Q$ and can be defined as follows: Let $\Omega=\{\omega: \R\to \R^d: \omega(\cdot) \,\,\mathrm{continuous}\}$ stand for all $\R^d$-valued continuous functions on $\R$, which, equipped with the topology of uniform convergence on bounded intervals, 
 is a complete separable metric space. The Borel $\sigma$-field of $\Omega$, denoted by $\mathcal F$, is generated by $\{\omega(t): -\infty<t<\infty\}$. Then we
 can identify  $\Omega $ with $\Omega_0\otimes \mathbb R^d$ by  mapping  
$$
\Omega\ni \omega\leftrightarrow (\omega^\prime,a) \qquad\mbox{where}\qquad a=\omega(0),\,\,\omega^\prime(t)=\omega(t)-\omega(0).
$$
Thus, any probability measure $\mathbb Q \in \mathcal M_1(\Omega)$ on $\Omega$ can be viewed as a measure on $\Omega_0\otimes \mathbb R^d$ and  will then have marginals ${\mathbb Q}^{\ssup 1} \in \mathcal M_{1}(\Omega_0)$, ${\mathbb Q}^{\ssup 2} \in \mathcal M_{1}(\mathbb R^d)$, respectively. The first marginal ${\mathbb Q}^{\ssup 1}$ is just the distribution of the increments of a process that has  $\mathbb Q$ for its distribution on $\mathcal F$.} i.e., $\mathbb Q$ is a stationary process on $\Omega$ and $\beta$ the distribution of its increments.
\end{Theorem}
Note that the above criterion does not need any tightness and is formulated in terms of convergence of integrals of continuous  functions vanishing at infinity w.r.t. measures on the function space $\Omega_0$. However, since $\Omega_0$ is not even locally compact, there is no notion of usual vague convergence of measures on this space (determined by convergence of integrals w.r.t. continuous functions vanishing at infinity). Therefore a notion of {\it wea-gue} convergence 
on measures on $\Omega_0\otimes \mathbb R^d$ was formulated in \cite[Section 2.1]{MV:18} which is conceptually important for the proof of the above result.

\section{Stochastic homogenization in the random conductance model}\label{sec-homogenization}

Part of the subsequent considerations may be viewed as a special instance of what has just been discussed before.

\vspace{.1cm}
The notion of a \emph{corrector} plays an important r\^ole in the context of \emph{stochastic homogenization} of a random media. We will describe the setup and how null-homology comes into play for a particular instance of a random walk in random environment (RWRE) in the reversible setup, known as  the \emph{random conductance model}. Let 
$$ E_d\ =\ \big\{(x,y)\colon |x-y|=1, \, x,y \in \Z^{d}\big\} $$ 
be the set of nearest neighbor bonds in $\Z^{d}$ and $\Omega= [a,b]^{E_d}$ for any two fixed numbers $0<a<b$. We assume that $\Omega$ is equipped with the product $\sigma$-field $\mathcal B$ and carries a probability measure $\Prob$. For simplicity, we also assume that the canonical coordinates are i.i.d. variables under $\Prob$. Note that any $x\in\Z^{d}$ acts on $(\Omega,\mathcal B,\Prob)$ as a $\Prob$-preserving and ergodic transformation $\tau_{x}$, defined as the canonical translation
$$ \Omega\ \ni\ \omega(\cdot)\,\mapsto\,\omega(x+\cdot). $$
For any $\omega\in \Omega$ and an associate family of probability measures $(\Prob_{\omega,x})_{x\in\Z^{d}}$, let then $(S_{n})_{n\geq 0}$ be the RWRE on $\Z^{d}$ which is defined as a Markov chain such that $\Prob_{\omega,x}(S_{0}=x)=1$ and transition probabilities are given by
\begin{equation}\label{pi}
\begin{aligned}
&\Prob_{\omega,x}(S_{n+1}=y+e|S_{n}=y)\ =\ \pi_\omega(y,y+e)\\
&\hspace{.3cm}:=\ \frac {\omega((y,y+e))}{\sum_{|e^\prime|=1} \omega((y,y+e^\prime))}\ 
=\ \pi_{\tau_{y}\omega}(0,e)
\end{aligned}
\end{equation}
for any $e$ with $|e|=1$ and $x\in\Z^{d}$. Furthermore, the sequence
$$ M_{n}\,:=\,\tau_{S_{n}}\omega\quad\text{for }n\ge 0, $$
with initial state $\omega$ and taking values in the ``environment space" $\Omega$, is also a Markov chain which drives the $S_{n}$, i.e.~$(M_{n},S_{n})_{n\ge 0}$ constitutes a MRW. It has transition kernel $P$ defined by
$$ (Pf)(\omega)\ :=\ \sum_{|e|=1} \pi_\omega(0,e) f(\tau_{e} \omega) $$
for all measurable, bounded $f$, is called the \emph{environment seen from the moving particle}, or just the \emph{environmental process}, and is particularly useful in the following scenario: Suppose there exists a probability density $\phi\in L^1(\Prob)$ (i.e. $\phi\geq 0$ and $\int \phi\,\dd\Prob=1$) such that $\Q=\phi\,\Prob$ is $P$-invariant, i.e.
\begin{gather}\label{invariance}
\langle Pf, \phi\rangle_{L^{2}(\Prob)}=\langle  f, \phi\rangle_{L^{2}(\Prob)},
\intertext{for all bounded and measurable $f$ or, equivalently,}
L^\star\phi\,=\,0,\quad L^{\star}\text{ the dual of }L\,=\,\mathrm{Id}-P.\nonumber
\end{gather} 
It can be shown, see \cite{PapVara:81,Kozlov:85,KipVara:86} and also \cite[Theorem 1.2]{BolSzn:02}, that such an invariant density $\phi$, if it exists, is necessarily unique. Moreover, $\Prob$ and $\Q$ are then  equivalent measures and $(M_{n})_{n\ge 0}$ as well as $(M_{n},X_{n})_{n\ge 1}$ (by Markov-modulation) are ergodic processes in equilibrium (under initial law $\Q$), where as usual $X_{n}=S_{n}-S_{n-1}$. 

\vspace{.1cm}
In the random conductance model with transition probabilities \eqref{pi}, the invariant density $\phi$ can actually be found by solving the detailed balance equations (reversibility), viz.
$$ \phi(\omega)\,=\,\frac 1C \sum_{|e|=1} \omega((0,e)),\quad\text{where}\quad C\,=\,\int\sum_{|e|=1}\omega((0,e))\ \Prob(\dd\omega). $$ 
Reversibility also provides that $P$ is self-adjoint on $L^{2}(\Q)$, that is 
\begin{equation}\label{reversible}
\langle f, P g\rangle_{L^{2}(\Q)}\ =\ \langle P f,g\rangle_{L^{2}(\Q)}
\end{equation}
for all bounded and measurable functions $f,g$.

\vspace{.1cm}
Returning to the RWRE $(S_{n})_{n\ge 0}$ under $\Prob_{\omega,0}$, the ergodicity of the Markov-modulated sequence $(M_{n},X_{n})_{n\ge 0}$ under $\Q$ easily provides that $S_{n}/n \to 0$ $\Prob_{\omega,0}$-a.s.~for $\Prob$-almost all $\omega$. To see this, let 
$$ \dd(\omega,x)\ =\ \Erw_{\omega,x}X_{1}\ =\ \sum_{|e|=1}e\pi_{\omega}(x,x+e)\ =\ \dd(\tau_{x}\omega,0) $$ 
denote the \emph{local drift at $x$} under $\Prob_{\omega,0}$. As $(S_{n})_{n\ge 0}$ has stationary ergodic increments under initial law $\Q$,  Birkhoff's ergodic theorem implies
$$ \frac{S_{n}}{n}\ \xrightarrow{n\to\infty}\ \int\Erw_{\omega,0}X_{1}\ \Q(\dd\omega)\ =\ \int\dd(\omega,0)\ \Q(\dd\omega)\ =\ 0 $$
for $\Q$-almost all and thus $\Prob$-almost all $\omega$ (as $\Prob,\Q$ are equivalent), the right-hand side being 0 by reversibility (recall \eqref{reversible}) and the definition of $\Q$. 

\vspace{.1cm}
As will be explained next, stochastic homogenization comes into play, and leads to the notion of a \emph{corrector}, when turning to the derivation of an almost sure CLT (or an invariance principle) for the law of $S_{n}$ under the quenched measure $\Prob_{\omega,0}$. Note that $Z_{n}= S_{n}-S_{0}-\sum_{j=0}^{n-1}\dd(\omega,S_{j})$, $n\ge 0$, is a $\Prob_{\omega,0}$-martingale with bounded increments (uniformly in $\omega$). Moreover, the local drift $\dd$ is bounded and thus particularly $\in L^{2}(\Prob)$. For any fixed $\eps>0$, let $g_{\eps}=\sum_{n\ge 1} \frac{P^{n-1}\dd}{(1+\eps)^{n}}$ be the $L^{2}(\Prob)$-solution to the perturbed Poisson equation 
$\big((1+\eps)\textrm{Id}-P\big)g_{\eps}=\dd $
that was also mentioned at the end of Subsection \ref{subsec-3.3}. Putting 
$$ G_{\eps}(\omega,e)\,:=\,(\nabla_e g_{\eps})(\omega)=g_{\eps}(\tau_e\omega)-g_{\eps}(\omega) $$ 
for any $e$ with $|e|=1$, Kipnis and Varadhan \cite[Theorem 1.3]{KipVara:86} showed that
$$ G_{\eps}(\cdot,e)\circ \tau_{x}\ \xrightarrow{L^{2}(\Prob)}\ G(\cdot,e)\circ\tau_{x}\quad\text{as }\eps\downarrow 0 $$ 
for any $x\in\R^{d}$, where $G$ is a (divergence free) 
\emph{gradient field}, i.e., it satisfies the \emph{closed loop condition}
\begin{equation}\label{closed loop condition}
\sum_{j=0}^{n-1}G(\tau_{s_{j}}\omega, {s_{j+1}-s_{j}})\ =\ 0\quad\Prob\text{-a.s.}
\end {equation}
for any closed path $s_{0}\to s_{1}\to\dots\to s_{n}=s_{0}$ in $\R^{d}$. The last property  allows us to define the \emph{corrector} corresponding to $G$ as 
\begin{equation}\label{corrector}
V_{G}(\omega,x)\,:=\,\sum_{j=0}^{n-1} G(\tau_{s_{j}}\omega, {s_{j+1}-s_{j}})
\end{equation}
along any path $0\to s_{1}\to\ldots\to s_{n-1}\to s_{n}=x$, the particular choice of the path being irrelevant because of \eqref{closed loop condition}. It also follows that $V_{G}$ has stationary and $L^{2}$-bounded gradient in the sense that
\begin{gather*}
V_{G}(\omega,y)-V_{G}(\omega,x)\ =\ V_{G}(\tau_{x}\omega,y-x)\quad\text{for all }x,y\in\Z^{d}
\shortintertext{and}
\sup_{x\in\Z^{d}}\|V_{G}(\cdot,x+e)-V_{G}(\cdot,x)\|_{L^{2}(\Prob)}\,<\,C,
\end{gather*}
respectively. Furthermore, fixing any $\omega\in\Omega$, the sequence
$$ (M_{n},V_{G}(\omega,S_{n}))_{n\ge 0} $$
forms a MRW under $\Prob_{\omega,0}$ whose driving chain is irreducible on the discrete state space $\Omega_{\omega}=\{\omega(x+\cdot):x\in\R^{d}\}$. The Markov-additive structure can be assessed by using \eqref{corrector}, which provides
$$ V_{G}(\omega,S_{n})\ =\ \sum_{j=1}^{n}G(\tau_{S_{j-1}}\omega,X_{j})\ =\ \sum_{j=1}^{n}G(M_{j-1},X_{j}) $$
for each $n\ge 0$, and by \eqref{closed loop condition} it also shows validity of the closed-loop condition. Hence, by invoking Theorem \ref{thm-main 4}, we infer that
\begin{equation}\label{eq:corrector NH}
V_{G}(\omega,S_{n})\ =\ \xi(M_{n})-\xi(M_{0})\ =\ \xi(\tau_{S_{n}}\omega)-\xi(\omega)\quad\Prob_{\omega,0}\text{-a.s.}
\end{equation}
for all $n\ge 0$ and a function $\xi:\Omega_{\omega}\to\R^{d}$. But this being true for each $\omega$, the function $\xi$ can be defined on the whole set $\Omega$ (in a measurable way) giving that \eqref{eq:corrector NH} holds $\Q$-a.s. In other words, the contractor is strict-sense null-homologous and thus tight, and it has a stationary version under $\Q$, namely $\xi(\omega)+V_{G}(\omega,S_{n})$ for $n\ge 0$. Although pointed out by Gloria \cite[p.~4]{Gloria:14} that stationary correctors in $L^{2}$ do not exist in dimension $d=1$ and $d=2$, this does not contradict our assertion. It rather implies that $\xi(\omega)$ cannot be square-integrable under $\Q$. In dimension $d\ge 3$, $L^{2}$-stationary correctors may exist under some additional conditions. Also the tightness is known in that case, see \cite{GloOtto:15} and \cite{ArKuuMou:19}.

\vspace{.1cm}
Now use that $x\mapsto V_{G}(\omega,x)+x$ is harmonic with respect to the transition probabilities \eqref{pi} for $\Prob$-almost all $\omega$ to infer that $(S_{n}+V_{G}(\cdot,S_{n}))_{n\ge 0}$ is a martingale with respect to $\Prob_{\omega,0}$. The corrector $V_{G}$ therefore expresses the ``distance" (or the \emph{deformation}) of the martingale from the random walk $(S_{n})_{n\ge 0}$ itself. The tightness ensures that the contribution of this deformation grows at most sub-linearly at large distances (i.e. $\sup_{|x|\leq n} n^{-1} V_{G}(x,\cdot)\xrightarrow{n\to\infty} 0$ a.s.) whence, by the martingale CLT, the laws $\Prob_{\omega,0}(S_{n}/\sqrt{n}\in \cdot)$ converge weakly to a Gaussian law for almost every $\omega$, see \cite{SidSzn:04,BerBis:07,MatPiat:07} for a detailed recount of the substantial progress made in this direction.

\vspace{.3cm}
\noindent{\bf Acknowledgment.}  The second author would like to thank S.R. S. Varadhan for many valuable discussions. Both authors would like to thank Sabine Jansen for pointing out the reference \cite[Theorem 3.1, Theorem 3.2]{AGL:01} (where a continuum version of the coboundary theorem has been proved using the ideas of Schmidt \cite{Schmidt:77}) and Karl Petersen for pointing out the references \cite{AR:19,P:73}. 


\bibliographystyle{abbrv}
\bibliography{StoPro}

\def\cprime{$'$}
\begin{thebibliography}{10}

\bibitem{AarWeiss:00}
J.~Aaronson and B.~Weiss.
\newblock Remarks on the tightness of cocycles.
\newblock {\em Colloq. Math.}, 84/85(part 2):363--376, 2000.
\newblock Dedicated to the memory of Anzelm Iwanik.


\bibitem{AR:19}
T. Adams and J. Rosenblatt.
\newblock{\em Existence and Non-existence of Solutions to the Coboundary Equation for Measure Preserving Systems.}
\newblock{\em Preprint.} arXiv: 1902.09045 (2019)



\bibitem{Alsmeyer:00}
G.~Alsmeyer.
\newblock The ladder variables of a {M}arkov random walk.
\newblock {\em Probab.~Math.~Statist.}, 20(1):151--168, 2000.

\bibitem{Alsmeyer:01}
G.~Alsmeyer.
\newblock Recurrence theorems for {M}arkov random walks.
\newblock {\em Probab.~Math.~Statist.}, 21(1):123--134, 2001.

\bibitem{AlsBuck:17c}
G.~Alsmeyer and F.~Buckmann.
\newblock Fluctuation theory for {M}arkov random walks.
\newblock {\em J.~Theoret.~Probab.}, 31(4):2266--2342, 2018.

\bibitem{AlsBuck:19}
G.~Alsmeyer and F.~Buckmann.
\newblock An arcsine law for {M}arkov random walks.
\newblock {\em Stochastic Process.~Appl.}, 129(1):223--239, 2019.

\bibitem{ArKuuMou:19}
S.~Armstrong, T.~Kuusi, and J.-C. Mourrat.
\newblock {\em Quantitative stochastic homogenization and large-scale
  regularity}, volume 352 of {\em Grundlehren der Mathematischen Wissenschaften
  [Fundamental Principles of Mathematical Sciences]}.
\newblock Springer, Cham, 2019.


\bibitem{AGL:01}
M.Aizenman, I. Goldstein, and J. Lebowitz
\newblock{\em Bounded fluctuations and translation symmetry breaking in one-dimensional 
particle systems.}
\newblock{\em J. Statist. Phys.} 103 (2001), no. 3-4, 601–618.



\bibitem{Benda:98}
M.~Benda.
\newblock A central limit theorem for contractive stochastic dynamical systems.
\newblock {\em J. Appl. Probab.}, 35(1):200--205, 1998.

\bibitem{BerBis:07}
N.~Berger and M.~Biskup.
\newblock Quenched invariance principle for simple random walk on percolation
  clusters.
\newblock {\em Probab. Theory Related Fields}, 137(1-2):83--120, 2007.

\bibitem{BlancLeBLions:07}
X.~Blanc, C.~Le~Bris, and P.-L. Lions.
\newblock The energy of some microscopic stochastic lattices.
\newblock {\em Arch. Ration. Mech. Anal.}, 184(2):303--339, 2007.

\bibitem{BolSzn:02}
E.~Bolthausen and A.-S. Sznitman.
\newblock {\em Ten lectures on random media}, volume~32 of {\em DMV Seminar}.
\newblock Birkh\"{a}user Verlag, Basel, 2002.

\bibitem{Bradley:95}
R.~C. Bradley.
\newblock On a theorem of {K}. {S}chmidt.
\newblock {\em Statist. Probab. Lett.}, 24(1):9--12, 1995.

\bibitem{Bradley:96}
R.~C. Bradley.
\newblock A ``multiplicative coboundary'' theorem for some sequences of random
  matrices.
\newblock {\em J. Theoret. Probab.}, 9(3):659--678, 1996.

\bibitem{Bradley:97}
R.~C. Bradley.
\newblock A ``coboundary'' theorem for sums of random variables taking their
  values in a {B}anach space.
\newblock {\em Pacific J. Math.}, 178(2):201--224, 1997.

\bibitem{DerriennicLin:01}
Y.~Derriennic and M.~Lin.
\newblock The central limit theorem for {M}arkov chains with normal transition
  operators, started at a point.
\newblock {\em Probab. Theory Related Fields}, 119(4):508--528, 2001.

\bibitem{DerriennicLin:01b}
Y.~Derriennic and M.~Lin.
\newblock Fractional {P}oisson equations and ergodic theorems for fractional
  coboundaries.
\newblock {\em Israel J. Math.}, 123:93--130, 2001.

\bibitem{DerriennicLin:03}
Y.~Derriennic and M.~Lin.
\newblock The central limit theorem for {M}arkov chains started at a point.
\newblock {\em Probab. Theory Related Fields}, 125(1):73--76, 2003.

\bibitem{Gloria:14}
A.~Gloria.
\newblock When are increment-stationary random point sets stationary?
\newblock {\em Electron. Commun. Probab.}, 19:no. 30, 14 pages, 2014.

\bibitem{GloriaOtto:11}
A.~Gloria and F.~Otto.
\newblock An optimal variance estimate in stochastic homogenization of discrete
  elliptic equations.
\newblock {\em Ann. Probab.}, 39(3):779--856, 2011.

\bibitem{GloOtto:15}
A.~Gloria and F.~Otto.
\newblock The corrector in stochastic homogenization: optimal rates, stochastic
  integrability, and fluctuations, 2015.
\newblock Preprint available at {\tt http://arxiv.org/abs/1510.08290}.

\bibitem{GordinLifsic:78}
M.~I. Gordin and B.~A. Lif{\v s}ic.
\newblock Central limit theorem for stationary {M}arkov processes.
\newblock {\em Dokl. Akad. Nauk SSSR}, 239(4):766--767, 1978.

\bibitem{Helland:82}
I.~S. Helland.
\newblock Central limit theorems for martingales with discrete or continuous
  time.
\newblock {\em Scand. J. Statist.}, 9(2):79--94, 1982.

\bibitem{Herrndorf:83}
N.~Herrndorf.
\newblock Stationary strongly mixing sequences not satisfying the central limit
  theorem.
\newblock {\em Ann. Probab.}, 11(3):809--813, 1983.

\bibitem{KipVara:86}
C.~Kipnis and S.~R.~S. Varadhan.
\newblock Central limit theorem for additive functionals of reversible {M}arkov
  processes and applications to simple exclusions.
\newblock {\em Comm. Math. Phys.}, 104(1):1--19, 1986.

\bibitem{Kozlov:85}
S.~M. Kozlov.
\newblock The averaging method and walks in inhomogeneous environments.
\newblock {\em Uspekhi Mat. Nauk}, 40(2(242)):61--120, 238, 1985.

\bibitem{Lalley:86}
S.~P. Lalley.
\newblock Renewal theorem for a class of stationary sequences.
\newblock {\em Probab. Theory Relat. Fields}, 72(2):195--213, 1986.

\bibitem{Leonov:61}
V.~P. Leonov.
\newblock On the dispersion of time-dependent means of a stationary stochastic
  process.
\newblock {\em Theory Probab.\ Appl.}, 6:87--93, 1961.

\bibitem{MatPiat:07}
P.~Mathieu and A.~Piatnitski.
\newblock Quenched invariance principles for random walks on percolation
  clusters.
\newblock {\em Proc. R. Soc. Lond. Ser. A Math. Phys. Eng. Sci.},
  463(2085):2287--2307, 2007.

\bibitem{MaxwellWood:00}
M.~Maxwell and M.~Woodroofe.
\newblock Central limit theorems for additive functionals of {M}arkov chains.
\newblock {\em Ann. Probab.}, 28(2):713--724, 2000.

\bibitem{MooreSchm:80}
C.~C. Moore and K.~Schmidt.
\newblock Coboundaries and homomorphisms for nonsingular actions and a problem
  of {H}. {H}elson.
\newblock {\em Proc. London Math. Soc. (3)}, 40(3):443--475, 1980.

\bibitem{MV:18}
C. Mukherjee and S. R. R. Varadhan.
\newblock{\em Identification of the Polaron measure in strong coupling and the Pekar variational formula.}
\newblock{\em Ann. Probab.}, 48(5):2119-2144, 2020


\bibitem{PapVara:81}
G.~C. Papanicolaou and S.~R.~S. Varadhan.
\newblock Boundary value problems with rapidly oscillating random coefficients.
\newblock In {\em Random fields, {V}ol. {I}, {II} ({E}sztergom, 1979)},
  volume~27 of {\em Colloq. Math. Soc. J\'{a}nos Bolyai}, pages 835--873.
  North-Holland, Amsterdam-New York, 1981.


\bibitem{P:73}
K. Petersen.
\newblock{\em On a series of cosecants related to a problem in ergodic theory.}
\newblock{\em Compositio Mathematica.} 26:3, 313-317, 1973.


\bibitem{Schmidt:77}
K.~Schmidt.
\newblock {\em Cocycles on ergodic transformation groups}.
\newblock Macmillan Company of India, Ltd., Delhi, 1977.
\newblock Macmillan Lectures in Mathematics, Vol. 1.
Available online: \href{https://www.mat.univie.ac.at/~kschmidt/Publications/rigveda.pdf}
{https://www.mat.univie.ac.at/~kschmidt/Publications/rigveda.pdf}



\bibitem{Shurenkov:84}
V.~M. Shurenkov.
\newblock On the theory of {M}arkov renewal.
\newblock {\em Theory Probab.\ Appl.}, 29(2):247--265, 1984.

\bibitem{SidSzn:04}
V.~Sidoravicius and A.-S. Sznitman.
\newblock Quenched invariance principles for walks on clusters of percolation
  or among random conductances.
\newblock {\em Probab. Theory Related Fields}, 129(2):219--244, 2004.

\bibitem{Woodroofe:92}
M.~Woodroofe.
\newblock A central limit theorem for functions of a {M}arkov chain with
  applications to shifts.
\newblock {\em Stochastic Process.~Appl.}, 41(1):33--44, 1992.

\bibitem{WuWood:00}
W.~B. Wu and M.~Woodroofe.
\newblock A central limit theorem for iterated random functions.
\newblock {\em J. Appl. Probab.}, 37(3):748--755, 2000.

\end{thebibliography}

\end{document}